\documentclass[reqno,fleqn,11pt]{amsart}
\usepackage{amsmath}
\usepackage{amssymb}
\usepackage{amsfonts}

\usepackage{graphicx}
\usepackage{amsthm}
\usepackage{enumerate}
\usepackage[T1]{fontenc}

\newcommand{\G}{\mathcal{G}}

\newcommand{\C}{\mathbb{C}}
\newcommand{\f}{\rightarrow}

\newcommand{\K}{K\"{a}hler}

\newcommand{\lmb}{\lambda}

\newcommand{\SL}{\operatorname{SL}}
\newcommand{\id}{\operatorname{id}}

\newcommand{\Aut}{\operatorname{Aut}}

\newcommand{\rank}{\operatorname{rank}}

\newcommand{\Ric}{\operatorname{Ric}}

\newtheorem{thm}{Theorem}
\newtheorem{prop}{Proposition}
\newtheorem{lem}[prop]{Lemma}

\newtheorem{exmp}[prop]{Example}

\newcommand{\fdim}{\hspace*{\fill}$\Box$}

\begin{document}
\title{Uniqueness of balanced metrics on holomorphic  vector bundles}
\author[A. Loi, R. Mossa]{Andrea Loi, Roberto Mossa}
\address{Dipartimento di Matematica e Informatica, Universit\`{a} di Cagliari,
Via Ospedale 72, 09124 Cagliari, Italy}
\email{loi@unica.it;   roberto.mossa@gmail.com }
\thanks{Research partially supported by GNSAGA (INdAM) and MIUR of Italy}
\date{February 25, 2010}

\begin{abstract}
Let $E\rightarrow M$ be a holomorphic vector bundle
over a compact \K\ manifold $(M, \omega)$.
We prove that if  $E$ admits  a   $\omega$-balanced  metric (in X. Wang's terminology \cite{wang})   then it is unique. This result   together  with \cite{ghigi}
implies  the existence  and uniqueness of $\omega$-balanced
metrics of certain direct sums of irreducible homogeneous vector bundles over rational homogeneous varieties.
We finally apply our result to show the rigidity of $\omega$-balanced \K\ maps into Grassmannians.
\end{abstract}

\maketitle

\section{Introduction and statement of the main results}

Let $E\rightarrow M$ be a very ample  holomorphic vector bundle over a compact
\K\ manifold $(M ,\omega)$ and let
$\underline s=(s_1, \dots , s_N)$ be a basis
of $H^0(M, E)$,  the space of  global holomorphic sections of $E$.
Let $i_{\underline s}:M\rightarrow G(r, N)$, $r=\rank E$,  be  the
Kodaira map associated to  the basis $\underline s$ (see, e.g. \cite{gh}), namely the  holomorphic embedding
whose expression $i_{\underline s}:U\rightarrow G(r, N)$
 in a local frame
  $(\sigma_1, \dots ,\sigma_r):U\rightarrow E$  is given by:
\begin{equation}\label{eqS}
i_{\underline s}(x)= \begin{bmatrix}
S_{11} (x)& \ldots & S_{1r}(x) \\
\vdots &  & \vdots \\
S_{N1}(x) & \ldots & S_{Nr}(x)
\end{bmatrix}, \ x\in U,
\end{equation}
where  $s_j=\sum_{\alpha=1}^rS_{j \alpha} \sigma_\alpha$, $j=1, \dots ,N$. Here
the square bracket  denotes the equivalence class in $G(r, N)=M^*(r,N, \C)/GL(r,\C)$,  where $M^*(r,N, \C)$ is the set of $r \times N$ complex  matrices  of rank $r$.

Consider the flat metric  $h_0$ on the tautological bundle $\mathcal T \rightarrow G(r, N)$,
i.e.  $h_0(v,w)=w^*v$,  and the
 dual metric  $h_{Gr}=h_0^*$ on the quotient bundle $\mathcal Q=\mathcal T^*$.
 Hence,  we can endow
 $E=i_{\underline s}^*\mathcal Q$ with the hermitian   metric
\begin{equation}\label{defhs}
h_{\underline s}=i_{\underline s}^*h_{Gr}
\end{equation}
 and define a   $L^2$-product on $H^0(M,E)$
by the formula:
\begin{equation}\label{deflhs}
< \cdot, \cdot >_{{h_{\underline s}}} =\frac{1}{V(M)} \int_M h_{\underline s}(\cdot, \cdot) \frac{\omega^n}{n!},
\end{equation}
where $\omega^n=\omega\wedge\cdots\wedge\omega$ and $V(M) = \int_M \frac{\omega^n}{n!}$.

An hermitian  metric $h$ on a very ample holomorphic vector bundle
$E\rightarrow M$ over a compact \K\ manifold $(M ,\omega)$
 is called $\omega$-{\em balanced} if there exists  a basis $\underline s$
 of $H^0(M, E)$
such that $h=h_{\underline s}=i_{\underline s}^*h_{Gr}$ and
\begin{equation}\label{balancedcondition}
<s_j, s_k>_{{h_{\underline s}}}=\frac{r}{N} \delta_{jk}, j, k=1, \dots , N=\dim H^{0}(M, E).
\end{equation}
One  also says that  $\underline s$ is
a $\omega$-{\em balanced} basis  of $H^0(M, E)$ if  (\ref{balancedcondition}) is satisfied. Therefore a metric on $E$ is
$\omega$-balanced if it is the pull-back of  the canonical metric $h_{Gr}$
of ${\mathcal Q}\rightarrow G(r, N)$ via the Kodaira map associated to a
$\omega$-balanced basis $\underline s$ of $H^0(M, E)$.
We refer the reader to  \cite{catlin}
for  the description of  interesting conditions on  a metric $h$ of a
holomorphic vector bundle $E\rightarrow M$ for the existence of  a (not necesseraly balanced) basis
$\underline s$ of $H^0(M, E)$ such that  $h=h_{\underline s}=i_{\underline s}^*h_{Gr}$.

The concept of balanced metrics on complex vector bundles  was introduced by X. Wang \cite{wang} (see also \cite{wang2}) following S. Donaldson's ideas \cite{donaldson}.
It  can be also  defined
in the non-compact case and the study of balanced  metrics
 is a very fruitful area of research both from
mathematical and physical point of view
(see, e.g.,  \cite{cgr3}, \cite{cgr4}, \cite{culoi}, \cite{baldisk}, \cite{mebal},
\cite{regscal}  and  \cite{albergbal}).

In \cite{wang}  X. Wang proves that, under the assumption that the \K\ form $\omega$ is integral,
$E$ is Gieseker stable if and only if  $E\otimes L^k$
admits a unique $\omega$-balanced metric (for every $k$ sufficiently large), where $L\rightarrow M$ is a polarization
of $(M, \omega)$, i.e. $L$ is a  holomorphic line bundle
over  $M$ such that $c_1(L)=[\omega]_{dR}$.

On the other hand, in  Lemma 2.7 of  \cite{reza},  R. Seyyedali shows  that if a simple bundle $E$
(i.e. $\Aut (E)=\C^*\id_E$, where $\Aut (E)$ denotes the group   of invertible holomorphic bundle maps from $E$ in itself) admits a balanced metric
then it is unique.
In the following  theorem, which is the main result of the present paper,
we prove the unicity of  balanced metrics  for  {\em any }  vector bundle.

\begin{thm}\label{maintheor}
Let $E$ be a holomorphic vector bundle over a compact \K\ manifold
$(M, \omega)$. If  $E$ admits a  $\omega$-balanced metric  then it is unique.
\end{thm}

As an  application of Theorem \ref{maintheor} and L. Biliotti and A. Ghigi results \cite{ghigi}
we obtain the existence  and uniqueness of $\omega$-balanced
metrics over certain direct sum of  homogeneous vector bundles over rational homogeneous varieties:

\begin{thm}\label{maincorol}
Let $(M, \omega )$ be a rational homogeneous variety and $E_j \f M, j=1, \dots, m$, be  irreducible homogeneous vector bundles over $M$ with $\rank E_j=r_j$ and $\dim H^0(M,E_j)=N_j> 0$, $j=1, \dots, m$.
If $\frac{r_j}{N_j}=\frac{r_k}{N_k}$ for all  $j, k=1, \dots, m$, then the homogeneous vector bundle $E=\oplus _{j=1}^m E_j\rightarrow M$ admits a unique
homogeneous $\omega$-balanced metric.
\end{thm}
\proof
Since
$\rank E=\sum_{j=1}^mr_j$ and
 $\dim (H^0(M, E))=\sum_{j=1}^mN_j$,
 it is enough to prove the theorem for $m=2$.
In \cite{ghigi} it is proved that each $E_j$, as in the statement, is a very ample bundle and  admits a unique homogeneous
$\omega$-balanced metric
 ${\underline s}^j=(s^j_1, \dots, s^j_{N_j})$, $j=1, 2$.
Then, the assumption $\frac{r_1}{N_1}=\frac{r_1}{N_2}$,
readily implies that  the basis
\begin{equation}\label{beq}
\underline s=((s^1_1,0), \dots, (s_{N_1}^1,0),(0,s_1^2),\dots, (0,s_{N_2}^2))
\end{equation}
is a homogeneous $\omega$-balanced basis for $E_1 \oplus E_2$. Then $h_{\underline s}=i_{\underline s}"^*h_{Gr}$ is the desired
homogeneous balanced metric on
$E_1 \oplus E_2$ which is unique   by Theorem \ref{maintheor}.
\fdim

\vskip 0.3cm

The proof of  Theorem  \ref{maintheor}  is based on Wang's work   on
balanced metrics (see \cite{wang} or the next section) and on  moment map techniques
developed by C. Arezzo and  the first author in \cite{arlcomm}, where it is proved   the unicity of   {\em balanced metrics} in the sense of S. Donaldson \cite{donaldson}.
Wang's work is summarized in the next section where we prove Lemma \ref{mainlemma} which is fundamental
for the proof of  Theorem \ref{maintheor},
to whom Section \ref{proofmaintheor} is dedicated.
In the last section we prove the rigidity of $\omega$-balanced \K\ maps into Grassmannians.

\vskip 0.3cm

\noindent
{\bf Acknowledgements:} We wish to thank Prof. Alessandro Ghigi and Prof. Reza Seyyedali for various interesting and stimulating discussions.

\section {Balanced bases and moment map}
Let $E$ be  a very ample   holomorphic vector bundle over a compact
\K\ manifold $(M, \omega )$.
Let  $J_0$  be the (complex) structure of $E$,
denote by $E_c$ the smooth complex vector bundle
underlying $E$ and  write $E=(E_c, J_0)$.
Let  $N$ be  the complex dimension of $H^0(M, E)$ and  let   $\mathcal H$ be the (infinite dimensional) manifold consisting of  pairs
$(\underline s,J)$ where $\underline s=(s_1, \dots, s_N)$ is an $N$-uple of complex  linearly independent smooth sections of $E_c$,  $J$ is a complex structure of $E_c$
and each section $s_j$ is holomorphic   with respect to the complex structure $J$, i.e.
\[
d s_j \circ I_0 = J \circ d s_j, \ j=1, \dots , N,
\]
where $I_0$ denotes the (fixed) complex structure of $M$.

Given  an hermitian metric $h$ on $E$  we denote by
$U_h (E_c)$ the subgroup of
$GL(E_c)$ consisting of smooth invertible bundle maps $E_c\rightarrow E_c$ preserving
the hermitian metric $h$ and by  $SU(N)\subset U(N)$ the group of  $N\times N$
unitary matrices with positive determinant. These groups act in a natural way on
$\mathcal H$ as follows:
\[
\Psi \cdot (\underline s,J) =(\Psi \underline s, \Psi\cdot  J),\  \Psi \in U_h (E_c)
\]
\[
U \cdot ( \underline s,J) = (U \underline s, J), \ U \in SU(N),
\]
were $\Psi \underline s = (\Psi s_1,\dots, \Psi s_N)$,
$\Psi\cdot J= \Psi J \Psi^{-1}$ and $U\underline s =
(U s_1,\dots, U s_N).$

Since these actions  commute  they induce a well-defined action of the group $\mathcal G_h = U_h(E_c) \times SU(N)$ on $\mathcal H$.
The   Lie algebra of $\mathcal G_h$
is   $GL (E_c) \oplus \mathfrak{su}(N)$ and its
complexification  $\mathcal G_h ^ \C=  GL (E_c) \times SL(N)$ naturally
 acts  on $\mathcal H$  by extending the action of $\mathcal G_h$.

\begin{thm}[Wang \cite{wang}]\label{wang}
The manifold $\mathcal H$ admits a \K\ form $\Omega$ invariant for the action of $\mathcal G_h$ whose moment map
$\mu_h: \mathcal H \f GL (E_c) \oplus \mathfrak{su}(N)$
is given by:
\begin{equation}\label{momentmapE}
\mu_h(\underline s,J)=(\ \sum_{j=1}^N h (\cdot, s_j)s_j,\  <s_j, s_k>_h - \frac{\sum_{j=1}^N
|s_j|^2_{h}}{N} \delta_{jk}\ ),
\end{equation}
where $|s_j|^2_{h}=<s_j, s_j>_h=\frac{1}{V(M)} \int_M h(\cdot, \cdot) \frac{\omega^n}{n!}$.
Consequently, a basis   $\underline s =(\underline s, J_0)$
of   $H^0(M, E)$ is balanced if and only if
$\mu_{h_{\underline s}}(\underline s, J_0)=(\id_E,0)$,
where  $h_{\underline s}$ is the metric of $E$ given by  (\ref{defhs}).
\end{thm}

A key ingredient in the proof of  Theorem \ref{maintheor} is  the following:
\begin{lem}\label{mainlemma}
Let  $\underline s=(\underline s, J_0)$ be a balanced basis of
$H^0(M, E)$ and let   $(\hat{\underline s}, \hat J)\in {\mathcal H}$
such that:
$\mu_{h_{\underline s}}(\hat {\underline s}, \hat J)=(\id_E,0)$
and $(\hat{\underline s}, \hat J)$ lies  in the same
$\mathcal G_{h_{\underline s}} ^ \C$-orbit of
$(\underline s, J_0)$.
Then
$(\hat{\underline s}, \hat J)$ lies  in the same
$\mathcal G_{h_{\underline s}}$-orbit of
$(\underline s, J_0)$, namely
 there  exists $(\Psi,U) \in \mathcal G_{h_{\underline s}}$ such that
$(\Psi,U)\cdot(\hat{\underline s}, \hat J) =(U\Psi\hat{\underline s}, \Psi\cdot \hat J)=(\underline s, J_0).$
\end{lem}
\proof
Since $a=(\id_E,0)\in GL (E_c) \oplus \mathfrak{su}(N)$ is (obviously) invariant by the coadjoint action  of
$\mathcal G_{h_{\underline s}}$, it is a standard fact in moment map's theory
 (see  \cite{hitchin} for a proof and also Proposition 3.1 in \cite{arlcomm}) that
$$\mu_{h_{\underline s}}^{-1}(a) \cap (\mathcal G_{h_{\underline s}}^ \C \cdot x)= \mathcal G_{h_{\underline s}} \cdot x, \ \forall x \in \mu_{h_{\underline s}} ^{-1}(a).$$
Then  the result follows by the assumptions and by Theorem \ref{wang}.
\fdim

\section{The proof of  Theorem \ref{maintheor}}\label{proofmaintheor}
Let  $E$ be  a very ample   holomorphic vector bundle over a compact
\K\ manifold $(M, \omega )$.
If $\underline s$ is {\em any}  basis of $H^0(M, E)$,
$F\in \Aut(E)$
and $U\in U(N)$,  then $i_{UF\underline s}=Ui_{F\underline s}=U i_{\underline s}$
where $UF\underline s =(UF s_1, \dots , UF s_N)$ and
it   follows easily  that  $h_{\underline s}=h_{UF\underline s}$.
Then the proof of  Theorem \ref{maintheor} will be a consequence of the following:
\begin{thm}\label{mainteorbal}
 If $s$ and $\tilde s$ are two balanced bases of $H^0(M, E)$ then there exist
a unitary matrix $U\in U(N)$ and $F\in \Aut (E)$
such that  $\underline{\tilde s}= U F \underline s$.
\end{thm}
\proof
Let  $h_{\underline s}$ and $h_{\tilde{\underline  s}}$
be the metric induced by $s$ and $\tilde s$ and
 $\Phi \in  GL (E_c)$ such that $\Phi ^*h_{\underline s}=h_{\tilde {\underline s}}$.
 We claim that
 \begin{equation}\label{equ1}
\sum_{j=1}^N h_{\underline s}(\cdot, \Phi \tilde { s}_j)\Phi \tilde { s}_j=\id_E
\end{equation}
and
\begin{equation}\label{equ2}
<\Phi \tilde{s}_j, \Phi \tilde {s}_k>_{h_{\underline s}}=\frac {r}{N} \delta_{jk}, \ j, k=1, \dots , N.
\end{equation}

Indeed
$$\id_E=\sum_{j=1}^Nh_{\underline {\tilde  s}}(\cdot ,  \tilde { s}_j)\tilde {s}_j=
\sum_{j=1}^N(\Phi^*h_{\underline s})(\cdot , \tilde { s}_j)\tilde {s}_j$$ and
if $\underline\sigma =(\sigma_1, \dots ,\sigma_r):U\rightarrow E$
is a local frame then, for all $\alpha =1, \dots ,r$,  one gets:
\begin{equation*}
\begin{split}
\sigma_{\alpha}=\Phi (\Phi^{-1}(\sigma_{\alpha}))&=
\Phi\left(\sum_{j=1}^Nh_{\tilde{\underline s}}(\Phi^{-1}(\sigma_{\alpha}) , \tilde { s}_j)\tilde {s}_j\right)\\
&=
\Phi\left(\sum_{j=1}^N(\Phi^*h_{\underline s})(\Phi^{-1}(\sigma_{\alpha}) , \tilde { s}_j)\tilde {s}_j\right)\\
&=\Phi\left(\sum_{j=1}^N h_{\underline s}(\sigma_{\alpha}, \Phi \tilde { s}_j)\tilde {s}_j\right)=
\sum_{j=1}^N h_{\underline s}(\sigma_{\alpha}, \Phi \tilde { s}_j)\Phi \tilde { s}_j,
\end{split}
\end{equation*}
where we have used the fact that
$\sum_{j=1}^N h_{\tilde {\underline s}} (\cdot, \tilde s_j)\tilde s_j=\id_E$,
and (\ref{equ1}) follows.

Moreover,
\begin{equation*}
\begin{split}
<\Phi \tilde s_j, \Phi \tilde s_k>_{h_{\underline s}}& = \frac{1}{V (M)} \int_M h_{\underline {s}}(\Phi \tilde s_j, \Phi \tilde s_k) \frac{\omega^n}{n!} \\
                                        &=\frac{1}{V (M)} \int_M (\Phi^*h_{\underline {s}})( \tilde s_j,  \tilde s_k) \frac{\omega^n}{n!}\\
                           & =\frac{1}{V (M)} \int_M h_{ {\tilde {\underline s}}}( \tilde s_j,  \tilde s_k) \frac{\omega^n}{n!} =\frac{r}{N} \delta_{jk}
\end{split}
\end{equation*}
and also  (\ref{equ2}) is proved.

It  follows by  (\ref{momentmapE}),  (\ref{equ1}) and (\ref{equ2})
that $({\underline s}, J_0)$ and $(\Phi \tilde {\underline s}, \Phi\cdot J_0)$ are in the same level set of the moment map
$\mu _{h_{\underline s}}$, namely
$$\mu _{h_{\underline s}}(\underline {s}, J_0)=\mu _{h_{\underline s}}(\Phi \underline {\tilde s}, \Phi\cdot J_0)=(\id_E, 0) .$$
Moreover,  since $\underline s$ and $\underline {\tilde s}$ are bases of the same vector space
$H^0(M, E)$ there exist a non zero constant  $\lambda$ and  $V \in \SL(N)$
such that $\lmb V \tilde {\underline s} = {\underline s}$.
Therefore
$$(\underline { s}, J_0)=(\lambda\Phi^{-1}, V)\cdot (\Phi\tilde {\underline s}, \Phi\cdot J_0)$$
and hence  $( {\underline s}, J_0)$ and $(\Phi  \tilde {\underline s}, \Phi J_0)$  are elements of ${\mathcal H}$  in the same $\G_{h_{\underline s}} ^\C$-orbit.
By  Lemma \ref{mainlemma}  there exists $(\Psi, U) \in \G_{h_{\underline s}} $ such that
\[
({\underline s},J_0)=(\Psi , U)\cdot(\Phi \tilde {\underline s}, \Phi\cdot J_0)=
(U\Psi\Phi \tilde {\underline s}, (\Psi\Phi)\cdot J_0).
\]
Consequently,   $F=\Psi\Phi:E_c\rightarrow E_c$ preserves the complex structure $J_0$, i.e.
$F\in \Aut (E)$ and ${\underline s}=UF \tilde {\underline s}$.
\fdim

\section{Rigidity of $\omega$-balanced \K\  maps into Grassmannians}
Let $(M, \omega)$ be a compact \K\ manifold.
A  holomorphic map $f:M\rightarrow G(r, N)$ is said to be $\omega$-balanced   if  there exist
a very ample  holomorphic vector bundle $E\rightarrow M$ and a balanced basis $\underline s$ of $H^0(M, E)$
such that  $f=i_{\underline s}$ (thus necessarily $f^*\mathcal Q=E$,  $r=\rank E$ and $N=\dim H^0(M, E)$).
A $\omega$-balanced map  $f:M\rightarrow G(r, N)$ is called a {\em \K\ } map
if $f^*\omega_{Gr}=\omega$, where $\omega_{Gr}$ is the standard \K\ form on $G(r, N)$, i.e.
$\Ric(h_{Gr})=\omega_{Gr}$.

\begin{exmp}\rm
Let $M={\C}P^1$ and $\omega_{\lambda} =\lambda \omega_{FS}$, where $\omega_{FS}$
is the  Fubini-Study \K\ form and $\lambda$ is  a positive real number.
Then, it not hard to see that the holomorphic map
\begin{equation}\label{balfac}
 f :{\C} P^1\rightarrow G(2, 4): [z_0, z_1]\mapsto \left[ \begin {array}{cc} z_0&0\\ \noalign{\medskip}0&z_0
\\ \noalign{\medskip}z_1&0\\ \noalign{\medskip}0&z_1\end {array} \right]
\end{equation}
is a  $\omega_{\lambda}$-balanced map for all  $\lambda$. Moreover, $f$ is \K\  when $\lambda=2$, i.e.
$f^*\omega_{Gr}=2\omega_{FS}$.
(In general it follows  by (\ref{deflhs}) and (\ref{balancedcondition}) that  if  $\underline s$
is a  $\omega$-balanced basis of $H^0(M, E)$ then $\underline s$ is still  $\lambda \omega$-balanced
for  $\lambda>0$).
\end{exmp}

Note that in the previous example $f^*\mathcal Q=O(1)\oplus O(1)$, where $O(1)$ is the hyperplane bundle on $\C P^1$,
$\Ric (h_{\underline s})=2\omega_{FS}$ (where $h_{\underline s}=f^*h_{Gr}$) and $f^*\omega_{Gr}=2\omega_{FS}$.
On the other hand, there exist  holomorphic maps $\tilde f=i_{\tilde {\underline s}}:\C P^1\rightarrow G(2, 4)$ (where $\tilde{\underline s}$
is a basis of $H^0\left(\C P^1,O\left(1\right)\oplus O\left(1 \right)\right)$)
satisfying these three  conditions  but for which  it cannot exist a  unitary transformation $U$ of $G(2, 4)$
such that $\tilde f=U f$ (cfr. \cite{chi}). An  example is given by:
\[
\tilde f :{\C} P^1\rightarrow G(2, 4): [z_0, z_1]\mapsto
\begin{bmatrix}
{z_0}^{2}&z_0\overline z_1 \frac{1}{2} \left(  \sqrt {3}-1 \right) \\
-z_0z_1\frac{1}{2} \left(  \sqrt {3}-1 \right) &|z_0|^2+\frac{1}{2} |z_1|^2\sqrt {3} \\
-z_0z_1 \frac{1}{2} \left( \sqrt {3}+1 \right) &-\frac{1}{2} |z_1|^2\\
{z_1}^{2}&\overline z_0z_1 \frac{1}{2} \left( 1- \sqrt {3}\right)
\end{bmatrix}.
\]

This phenomenon is due to the  fact that
the rigidity of   \K\ maps into   $G(r, N)$ with $r\geq 2$
does not in general hold true
 (see, e.g.,  \cite{calabi},  \cite{chi}, \cite{green},  \cite{peng1}, \cite{peng2}),
in contrast with the case $r=1$ where one has  the celebrated Calabi's rigidity theorem
for \K\ maps into projective spaces.

\vskip 0.3cm

On the other hand the following theorem, which is the main result of this section, shows the rigidity of $\omega$-balanced \K\ embedding.

\begin{thm}\label{teorrig}
Let  $E\rightarrow M$ be a very ample complex vector bundle
over a compact \K\ manifold $(M, \omega )$. Assume that $E$ admits a  $\omega$-balanced metric
$h$
such that $\Ric (h)=\omega$.
Then  there exists a unique (up to a unitary transformations of $G(r, N))$
$\omega$-balanced \K\ embedding
$f:M\rightarrow G(r, N)$ such that
$f^*Q=E$.
\end{thm}
\proof
Let $\underline s$ be a balanced basis
of $H^0(M, E)$ and let  $f=i_{\underline s}:M\rightarrow G(r, N)$
be the associated Kodaira's  map.
By Theorem \ref{maintheor} $f$ is the unique
 (up to a unitary transformations of $G(r, N)$)
$\omega$-balanced
embedding such that $f^*Q=E$.
So it remains to show that   $f^*\omega_{Gr}=\omega$.
Fix a local frame   $(\sigma_1, \dots ,\sigma_r):U\rightarrow E$.
In this local frame   $f:U\rightarrow G(r, N)$ is given by
(\ref{eqS}).
Then, the local expression of $h=h_{\underline s}$,
$\omega=\Ric (h)$
and $f^*\omega_{Gr}$
are given respectively by
$(S^*S)^{-1}$,
$-\frac{i}{2}
\partial\bar\partial\log\det (S^*S)^{-1}$
and
$\frac{i}{2}
\partial\bar\partial\log\det (S^*S)$.
\fdim

\end{document}